\def\draft{n}
\documentclass{amsart}
\usepackage{fullpage,epic,eepic,epsfig,amssymb}
\theoremstyle{plain}

\newtheorem{theorem}{Theorem}
\newtheorem{proposition}{Proposition}[section]
\newtheorem{lemma}[proposition]{Lemma}
\newtheorem{corollary}[proposition]{Corollary}

\theoremstyle{definition}
\newtheorem{definition}[proposition]{Definition}

\theoremstyle{remark}

\def\printname#1{
	\if\draft y
		\smash{\makebox[0pt]{\hspace{-0.5in}
			\raisebox{8pt}{\tt\tiny #1}}}
	\fi
}

\newcommand{\psdraw}[2]
         {\begin{array}{c} \hspace{-1.3mm}
	\raisebox{-4pt}{\epsfig{figure=draws/#1.eps,width=#2}}
	\hspace{-1.9mm}\end{array}}

\newlength{\standardunitlength}
\catcode`\@=11
\long\def\@makecaption#1#2{%
     \vskip 10pt

\setbox\@tempboxa\hbox{%\ifvoid\tinybox\else\box\tinybox\fi
       \small\sf{\bfcaptionfont #1. }\ignorespaces #2}%
     \ifdim \wd\@tempboxa >\captionwidth {%
         \rightskip=\@captionmargin\leftskip=\@captionmargin
         \unhbox\@tempboxa\par}%
       \else
         \hbox to\hsize{\hfil\box\@tempboxa\hfil}%
     \fi}
\font\bfcaptionfont=cmssbx10 scaled \magstephalf
\newdimen\@captionmargin\@captionmargin=2\parindent
\newdimen\captionwidth\captionwidth=\hsize
\catcode`\@=12

\def\lbl#1{\label{#1}\printname{#1}}

\def\eqdef{\overset{\text{def}}{=}}

\def\BN{\mathbb N}
\def\BZ{\mathbb Z}
\def\BQ{\mathbb Q}
\def\BR{\mathbb R}

\def\A{\mathcal A}
\def\B{\mathcal B}

\def\S{\Sigma}

\def\ihs{integral homology 3-sphere}
\def\qhs{rational homology 3-sphere}

\def\fti{finite type invariant}

\def\la{\langle}
\def\ra{\rangle}

\def\o1o{\underset{1}\ast}
\def\x1y{\underset{x1y}\ast}
\def\y1x{\underset{y1x}\ast}

\def\FY#1{\mathcal F^{\mathrm{Y}}_{#1}(M)}

\def\GY#1{\mathcal G^{\mathrm{Y}}_{#1}(M)}
\def\spanning{$H_1$-spanning}

\def\qbasis{$H_{1,\BQ}$-basis}

\def\e{\varepsilon}

\def\AS{\mathrm{AS}}
\def\IHX{\mathrm{IHX}}
\def\LOOP{\mathrm{LOOP}}

\def\OBRA{\mathrm{OBR}}
\def\CBRA{\mathrm{BR}}
\def\clover{clover} 
\def\ygraph{$\mathrm{Y}$-graph}

\def\Ao{{\mathcal A}^{\mathrm o}}
\def\Ac{{\mathcal A}}

\begin{document}

%%%%%%%%%%%%%%%%%%%%%%\inplude{page1}

\title[The mystery of the brane relation]{The mystery of the brane 
relation}

\author{Stavros Garoufalidis}
\address{School of Mathematics \\
          Georgia Institute of Technology \\
          Atlanta, GA 30332-0160, USA. }
\email{stavros@math.gatech.edu}

\thanks{The  author partially supported by NSF grant
        DMS-98-00703  and by an Israel-US BSF grant.\newline
        This and related preprints can also be obtained at
{\tt http://www.math.gatech.edu/$\sim$stavros } 
\newline
1991 {\em Mathematics Classification.} Primary 57N10. Secondary 57M25.
}

\dedicatory{Dedicated to the memory of M. Goussarov.}

\date{
This edition: June 6, 2000 \hspace{0.3cm} First edition: October 27, 1999.}

\begin{abstract}
The purpose of the present paper is to 
introduce and explore two surprises that arise when we apply a standard
procedure to study the number of 
\fti s of 3-manifolds introduced independently
by M. Goussarov and K. Habiro based on surgery on claspers, \ygraph s or
\clover s, \cite{Gu,Ha,GGP}.
One surprise is that the upper bounds depend on a bit more than a choice
of generators for $H_1$. A complementary surprise a curious
brane relation (in two flavors, open and closed)
which shows that the upper bounds are in a certain
sense independent of the choice of generators of $H_1$.
\end{abstract}

\maketitle

%\tableofcontents

%%%%%%%%%%%%%%%%% the text file

\section{Introduction}
\lbl{sec.intro}

\subsection{Motivation}
\lbl{sub.motivation}

It is well-known that starting from a {\em move} (often described in terms
of surgery on a link) on a set of knotted objects 
(such as knots, links, braids,
tangles, 3-manifolds, graphs), one can define a notion
of {\em \fti s}. The question of how many invariants are there in any degree
gets divided into two separate questions: one that provides upper bounds for 
the number of invariants, and one that provides lower bounds.
Traditionally, upper bounds are obtained by providing
a set of topological relations
among the moves, whereas lower bounds are obtained by constructing
(by quite different means) invariants. 
The purpose of the present paper is to 
introduce and explore two surprises that arise when we study the first question
regarding the \fti s of 3-manifolds due, independently,
to M. Goussarov and K. Habiro based on surgery on \clover s, 
see \cite{Gu,Ha}.
One surprise is that the upper bounds depend on a bit more than a choice
of generators for $H_1$. A complementary surprise a curious
brane relation which shows that the upper bounds are in a certain
sense independent of the choice of generators of $H_1$.

Let us briefly recall that Goussarov and Habiro
studied the notion of surgery $M_G$ of a 3-manifold $M$ along an (embedded)
{\em clasper} or  $\mathrm{Y}$-{\em graph} $G$, see \cite{Gu,Ha,GGP}. 
For a detailed definition of \clover s
 and their associated surgery see \cite{Gu,Ha} or 
\cite{GGP} the notation of which we follow here. 
This notion of surgery allows one to consider the abelian group $\FY {}$
freely generated by (isomorphism classes of) 3-manifolds obtained surgery on
$M$, and define
a decreasing filtration on it
where $\FY m$ is the subgroup generated
by 
$$[M,G]=\sum_{G' \subset \{G_1, \dots, G_n \}} (-1)^{|G'|} M_{G'}$$
for all \clover s $G$ in $M$ partitioned into $n \geq m$ blocks
$G_1, \dots, G_n$, where for $G' \subset \{G_1, \dots, G_n \}$,
$M_{G'}$ denotes the result of surgery on the graph 
$\cup_{i: G_i \in G'} G_i$, and where $|G'|$ denotes the cardinality of
the set $G'$.

Dually, and perhaps more naturally, this filtration allows us to call
a $\BZ$-valued function $\lambda$ on the set of 3-manifolds obtained
by surgery on \clover s a \fti \ of type $m$ iff $\lambda(\FY {m+1})=0$.
Thus, the question of how many \fti s of type $m$ are there translates into
a question about the structure of the (graded quotient) abelian groups
$\GY m\eqdef\FY m/\FY {m+1}$. For the case of $M=S^3$ (or any other \ihs ),
it is well-known that the topological calculus of \clover s developed 
indepedently by Goussarov and Habiro, implies the existence of upper bounds 
of $\GY {}$ in terms of an abelian group $\A (\phi)$ generated by (abstract) 
trivalent graphs, modulo the well known antisymmetry $\AS$ and $\IHX$ 
relations, see for instance \cite[Section 4]{GGP}. The case of arbitrary 
3-manifolds $M$ (needed for instance in \cite[Theorems 5,6,7]{GL})  
seems to be missing from the literature, even though the main tools are the 
same as in the case of $M=S^3$. There are, however, two surprises
in extending the above upper bound to all 
closed 3-manifolds, which are the main point of this paper: one is 
that the upper bound for $\GY {}$ is given in terms of a
finitely generated (in each degree)
abelian group $\Ao (b)$ defined below, where $b$ is a 
{\em \spanning } link i.e., an oriented
framed link in $M$ that generates (possibly with redundances)
$H_1(M,\BZ)$, see Theorem \ref{thm.1}. In other words, the generators of
$\Ao (b)$ {\em depend} on just a bit more than a
choice of generators for $H_1(M,\BZ)$, they depend on a choice of 1-cycles. 
The other surprise is the existence of a new relation in $\Ao (b)$, the 
{\em open brane} ($\OBRA$) and the {\em closed brane} ($\CBRA$) relation, 
which is also given in terms of a choice of embedded 2-cycles in $M$.

Of course the choice of $b$ is not unique, and the choice of cycles
in the $\OBRA$ relation is not unique, however the $\OBRA$ and $\CBRA$ 
relations
imply that any two such choices $b$ and $b'$ lead to rather canonical
isomorphisms between $\Ao(b)$ and $\Ao(b')$ as well as commutative
diagrams, see Theorem \ref{thm.1}.

If one is willing to work with rational coefficients, then the above
upper bound $\Ao(b)$ can be identified with an invariant $\A$-group
$\Ac(H(M))$ that depends only
on the cohomology ring $H^\ast(M,\BQ)$ of $M$, see Corollary \ref{cor.Q}
(although the map $\Ac(H(M))\to\GY {}$ still depends on a choice
of a \spanning \ link $b$).

As a final comment before the details, we should mention that for \fti s
of \ihs s, or for $\BQ$-valued \fti s of rational homology spheres
the above mentioned choices of 1-cycles and 2-cycles are invisible, which
partly explains why they were not discovered despite the success of 
constructing and studying theories of \fti s of rational homology spheres.

\subsection{Statement of the results}
\lbl{sub.results}

Throughout, by {\em graph} we mean 
we mean one with (symmetric) univalent and trivalent 
vertices, together with a choice of cyclic order on each trivalent vertex.
Note that graphs that contain struts, i.e., an interval with two univalent
vertices and no trivalent ones, will {\em not} be allowed here.
Univalent vertices of graphs will often be called {\em legs} or {\em leaves}. 
Given a set $X$, an $X$-{\em colored graph} is a graph $G$ together with
a function $c:\mathrm{Legs}(G)\to X$. 
This assignment can be extended
linearly to include graphs each univalent vertex of which is assigned
a nonzero formal linear combination of elements of $X$.

Let $\B(X)$ denote the abelian group spanned by $X$-colored graphs
modulo the well-known $\AS, \IHX$ and $\LOOP$ relations shown in Figure
\ref{relations}. $\B(X)$ is graded, by declaring the degree of a graph
to be the number of its trivalent vertices.

Notice that the group $\B(X)$ is closely related to a group that appears when
one studies \fti s of $X$-component links in $S^3$, with some notable 
differences: one is that we do not allow struts, another is that we do not 
grade by half the number of vertices, and the third is that we allow graphs
with no legs.
  
\begin{figure}[htpb]
$$ \psdraw{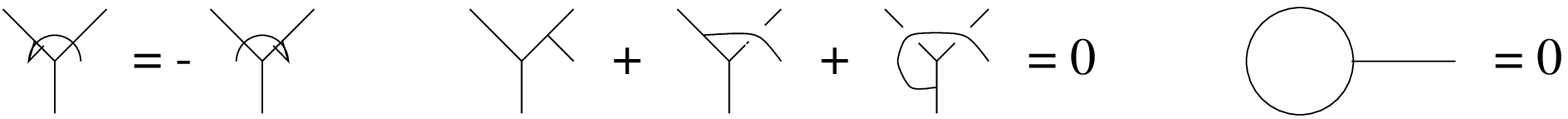}{5in} $$
\caption{The $\AS$, $\IHX$ (all trivalent vertices oriented counterclockwise),
 and $\LOOP$ relations. In the $\LOOP$ relation, the appearing loop
is an edge and not a leaf of the graph.}\lbl{relations}
\end{figure}

Given a \spanning \ link $b$, we now define two important 
relations on $\B(b)$. Let $\cdot: H_2(M,\BZ)\otimes 
H_1(M,\partial M,\BZ)\to\BZ$ be the 
intersection pairing. 

\begin{definition}
\lbl{def.cbra}
Fix a closed surface $\S$ in $M$.
Let $(G,\ast)$ be $b$-colored graph, which contains a special
leg colored by the special symbol
$\ast$ (disjoint from the alphabet $b$).
Let 
$$\langle G, \S \rangle :=\sum_{l} [\S]\cdot [c_l] G_l \in \B (b)$$ 
where the summation is over all legs of $G$ except $\ast$
and 
where $G_l$ is the result of gluing 
the $\ast$-leg of $G$ to a  $c_l$-colored leg of 
$G$, as shown in the following example 
$$\left\la \psdraw{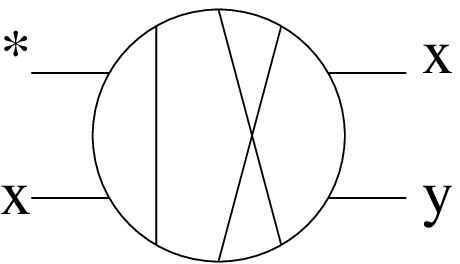}{0.8in}, \S \right\ra=
[\S]\cdot[x]
\psdraw{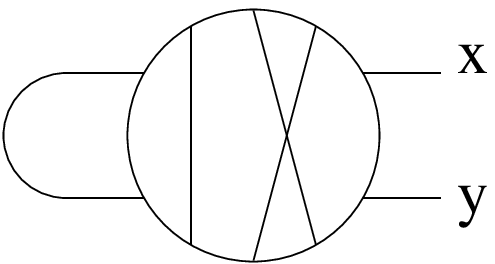}{0.8in} + [\S]\cdot[x]
\psdraw{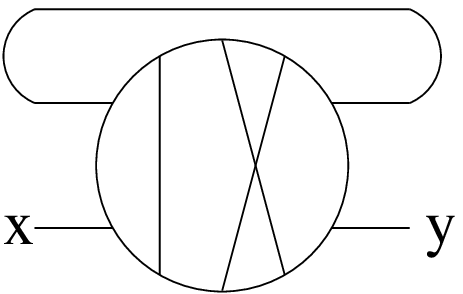}{0.8in} + [\S]\cdot[y]\psdraw{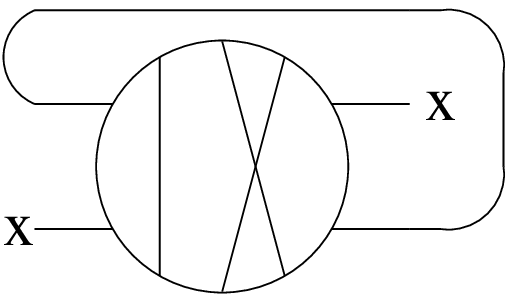}{0.8in}.$$
By convention, the summation over the empty set equals to zero. 
The $\CBRA$ ({\em closed brane}) relation\footnote{which does not seem to be
related in any meaningful way to the wonderful (mem)branes of string theory.}
is the subgroup of $\B(b)$ generated by $\la G, \S \ra=0$ for all surfaces
$\S$, or really, only a generating set for $H_2(M,\BZ)$
and all graphs $(G,\ast)$ as above. Let $\Ac(b)=\B(b)/(\CBRA)$.
\end{definition}

\begin{definition}
\lbl{def.obra}
Fix a $b$-colored graph that contains a distinguished leg $\ast$ colored
by a nullhomologous label $c_0$ which bounds a surface $\S_0$ in $M$.
Let 
$$ \la G, \S_0 \ra:=G+ \sum_l [\Sigma_0]\cdot [c_l] G_l \in \B(b)$$
where the summation is over all legs of $G$ except $\ast$ and 
where $G_l$ is the result of gluing the $\ast$-leg of $G$ to a $c_l$-colored
leg of $G$. The $\OBRA$ ({\em open brane}) relation is the subgroup of $\B(b)$
generated by $\la G, \S_0 \ra=0$ for all graphs $G$ as above and all surfaces
$\S_0$. Note that the $\OBRA$ subgroup of $\B(b)$ includes the $\CBRA$
subgroup if we assume that one of the components of $b$ is the boundary
of an embedded disk disjoint from the rest of the components of $b$. Let 
$\Ao(b)=\Ac(b)/(\OBRA)$.
\end{definition}

\begin{theorem}
\lbl{thm.1}
(i) For every \spanning \ link $b$ in a manifold $M$, 
there is a group homomorphism
$$W_{M,b}: \Ao (b) \to \GY {}$$
which is onto, once tensored with $\BZ[1/2]$. \newline
(ii) For every two \spanning \ links $b$ and $b'$ in $M$, there are 
isomorphisms 
$W_{M,b,b'}: \Ao (b) \to \Ao (b')$ such that over $\BZ[1/2]$:
\begin{equation}
\lbl{eq.1}
 W_{M,b}=W_{M,b'} \circ W_{M,b,b'}.
\end{equation}
\end{theorem}

\subsection{The size of $\Ao (b)$}
\lbl{sub.size}
It is natural to ask how big is the (finitely generated in each 
degree) abelian
group $\Ao (b)$ which bounds from above $\GY {}$. 

\begin{corollary}
\lbl{cor.1}
(i) If $H_1(M,\BZ)$ is torsion-free and $b$ is a basis of $H_1$, then
$$
\Ac(b) \cong \Ao (b).$$
(ii) If $b$ is \qbasis \ and $b'$ is \spanning \ then
$$ \Ac(b) \cong_{\BQ} \Ao(b'). $$
(iii)
If %the intersection form $H_2(M,\BQ) \otimes H_1(M,\partial M,\BQ) \to \BQ$
%of a 3-manifold $M$ vanishes, 
$H_1(M,\partial M,\BQ)=0$, then for every \spanning \ $b$ we have
$$\B(b)\cong \Ac(b).$$
(iv)
In particular, for $M$ a \qhs , we have that
$$
\Ao(b)\cong_{\BQ} \Ac(b)\cong_{\BQ} \A(\phi).
$$
(v)
For $M$ a homology-cylinder (i.e., a manifold with the same integer homology 
as that of $\S\times I$ for a surface $\S$ with one boundary component)
and a \qbasis \ $b$, we have that
$$
\B(b)\cong_{\BQ} \Ac(b)\cong_{\BQ} \Ao(b).
$$
\end{corollary}

If we are willing to work with rational coefficients, then one can
define in an invariant way a group of graphs, that depends only on
the cohomology ring $H^\ast(M,\BQ)$ as follows:
$\Ac(H(M))$ is generated by graphs colored by nonzero elements of
$H_1(M,\BQ)$, modulo the $\AS$, $\IHX$, $\LOOP$ and $\CBRA$ relations.

\begin{corollary}
\lbl{cor.Q}
For every manifold $M$, there is a map
$$\Ac(H(M))\to \GY{},$$ 
onto over $\BQ$.
\end{corollary}

For manifolds $M$ with $b_1(M)=0$, i.e., for \qhs s, we show a promised 
{\em universal property} of the LMO invariant restricted to the set of 
rational homology spheres \cite{LMO}, or of its cousin, the Aarhus integral 
\cite{A}:

\begin{theorem}
\lbl{thm.2}
The LMO invariant is the universal $\BQ$-valued  \fti \ of
rational homology spheres. In particular, for $M$ a \qhs \ and $b$ \spanning ,
we have $\Ao (b) \cong_{\BQ} \Ac (b) \cong_{\BQ} \A (\phi) \cong_{\BQ}
\GY {}$. 
\end{theorem}

With regards to the size of $\A (\phi)$, it is well-known that 
Lie algebras and their representation theory provides lower bounds
for the abelian groups $\A (\phi)$. 
In the case of manifolds $M$ with positive betti number,
we do not know yet 
of lower bounds for $\Ao (b)$. The little we know at present is the following:

\begin{corollary}
\lbl{cor.3}
Let $b$ be a \spanning \ link in a closed manifold $M$ and $G$ be a 
graph colored by a sublink $b'$ of $b$. Assume that
$G$ has an internal edge, that is an edge between two trivalent 
vertices of $G$.
If $b'$ is not \spanning \ (over $\BQ$), then $G=0 \in \Ao(b)$. \newline
In particular, if $b_1(M) > 0$, then every graph without legs vanishes
in $\Ao (b)$.
\end{corollary}

\begin{corollary}
\lbl{cor.levine}
Let $b$ be a \spanning \ link in a closed manifold $M$ and
$G$ be a graph whose $r+1$
legs are colored by $x, y_1,\dots,
y_r$ so that $x$ is primitive and linearly independent from 
$\{y_1,\dots,y_r\}$. For $k=0,\dots , r$, let $G^{(k)}$ denote the sum of all 
ways of replacing $k$ many $y_i$ by $x$. Assume that $G$ contains an internal
edge. Then $G^{(k)}=0$ in $\Ao (b)$ for all $k$.
\end{corollary}

We caution the reader that the above corollary by no means implies that
$\Ao (b)$ is zero dimensional for manifolds $M$ with positive betti number,
since for example, for manifolds with positive betti number, the coefficients
of the Alexander polynomial (of the maximal torsion-free abelian cover)
are \fti s in our sense.

\section{Proofs}
The proofs of the theorem and its corollaries involve
standard alternations of the {\em topological calculus of 
\clover s} as appears in \cite{Gu,Ha}; the uninitiated reader
may also look at \cite[Section 3]{GGP}.

Before we prove the theorem, it will be important to state 
some lemmas
the proof of which follows by applying to the topological calculus of \clover s
elementary alternations, see for example \cite[Section 4.1]{GGP}:

\begin{lemma}\cite{Ha,Gu}{\bf (Cutting a Leaf)}
\lbl{lem.cut}
Let $G$ be a \clover \ of degree $m$ in a manifold $M$ 
and $L$ be a leaf of $G$. 
An arc $a$ starting in the external vertex incident to $L$ and 
ending in other point of $L$, splits $L$ into two arcs $L'$ and $L''$.
Denote by $G'$ and $G''$ the graphs obtained from $G$ by replacing 
the leaf $L$ with $L'\cup a$ and $L''\cup a$ respectively, see Figure
\ref{split}.
Then $[M,G]=[M,G']+[M,G'']$ in $\GY {m}$.
\end{lemma}

\begin{figure}[htpb]
$$ \psdraw{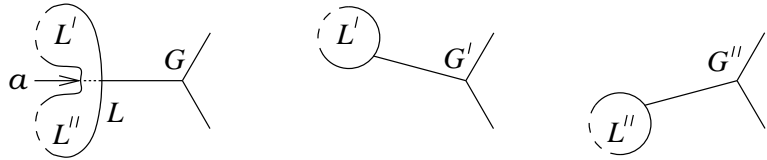}{4.0in} $$
\caption{Splitting a leaf.}\lbl{split}
\end{figure}

\begin{lemma}\cite{Gu,Ha}{\bf (Sliding an Edge)}
\lbl{lem.slide}
Let $G$ be a \clover \ of degree $m$ in a manifold $M$, and let
$G'$ be obtained from $G$ by sliding an edge of $G$ along a tube in $M$.
Then $[M,G]=[M,G']$ in $\GY m$.
\end{lemma}

\begin{lemma}\cite{Gu,Ha}
\lbl{lem.relations}
For all $\e=\pm 1$, we have the following identities in $\GY {}$: 
\begin{equation*}
2 [M, \psdraw{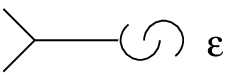}{0.8in}]  =  0  \quad\quad \text{ and } \quad\quad
[M,\psdraw{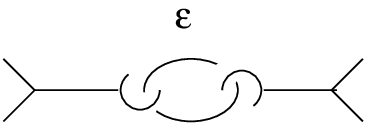}{0.9in}]  =  -\e [M,\psdraw{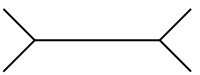}{0.5in}].
\end{equation*}
\end{lemma} 

\begin{lemma}
\lbl{lem.obra}
Let $G$ be a \clover \ with $r+1$ leaves
$l_i$ for $i=0, \dots, r$ in a manifold $M$.
Assume that $l_0$ bounds an embedded surface $\Sigma_0$ in $M$.
Then
$$ G+ \sum_{i=1}^r [\Sigma_0]\cdot[l_i] G_i=0 \in \GY {}$$
where $G_i$ is the result of gluing the $0$-th leg of $G$
to its $i$-th leg.
\end{lemma}

\begin{proof}
Consider a graph $G$ and a surface $\Sigma_0$ as above.
$\Sigma_0$ can be thought of as an embedded disk
with bands. We can assume that $G$ is disjoint 
from the (interiors of the bands) of $\Sigma_0$ and thus $G$ intersects
the (interior of) $\Sigma_0$ only in the embedded disk.
Cut each band along arcs (in the normal direction to the core
of the band) using the Cutting and Sliding Lemmas \ref{lem.cut} and 
\ref{lem.slide} as shown
\begin{eqnarray*}
[M,\psdraw{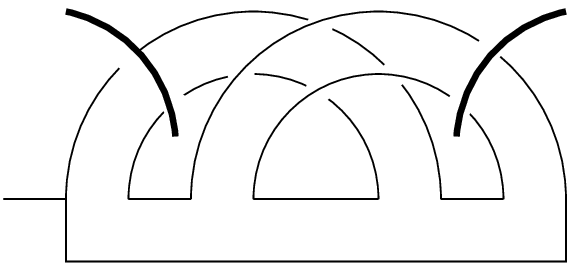}{1.1in}]  & = &  [M,\psdraw{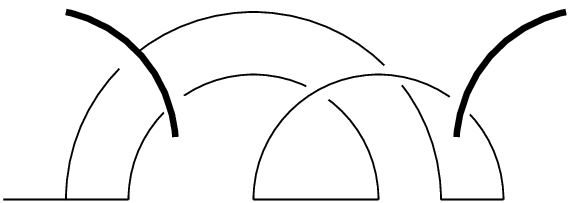}{1.1in}]
+ [M,\psdraw{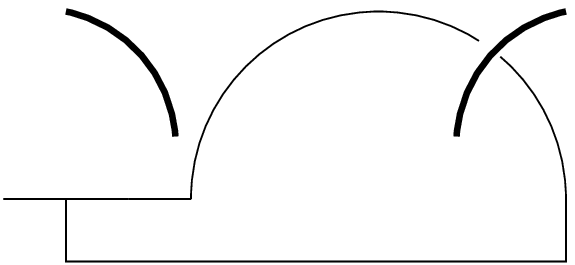}{1.1in}] \\
& = & [M,\psdraw{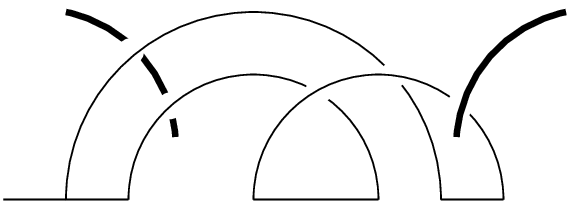}{1.0in}]  +  [M,\psdraw{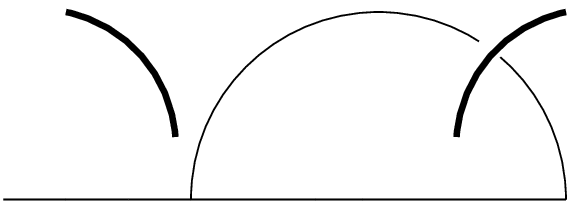}{1.0in}]
+ [M,\psdraw{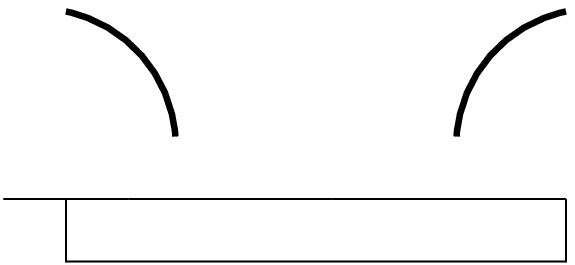}{1.0in}] \\
& = & [M,\psdraw{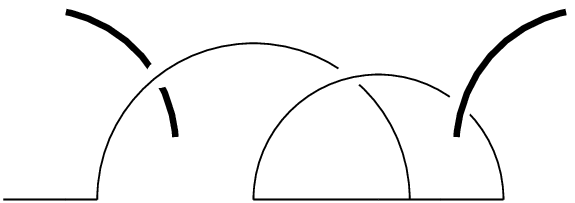}{1.0in}]  +  [M,\psdraw{Cutsurface5}{1.0in}]
+ [M,\psdraw{Cutsurface6}{1.0in}] \\
& = & -[M,\psdraw{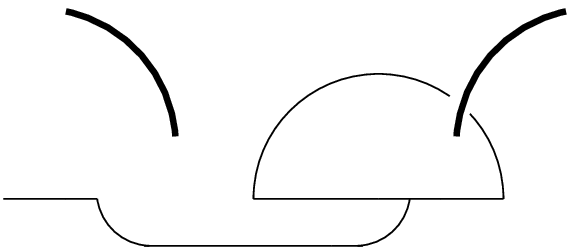}{1.0in}]  +  [M,\psdraw{Cutsurface5}{1.0in}]
+ [M,\psdraw{Cutsurface6}{1.0in}] \\
& = & [M,\psdraw{Cutsurface6}{1.0in}]
\end{eqnarray*}
(where $\S_0$ is a surface of genus $1$, and the solid arcs represent 
arbitrary tubes in the 3-manifold). The above calculation reduces to the
case of a surface $\Sigma_0$ hawith no bands, i.e.,  a disk. Using the 
Cutting and Sliding Lemmas \ref{lem.cut} and \ref{lem.relations} once again,
we may assume that the leaf $l_0$ of $G$ is zero-framed and that the disk
$\Sigma_0$ intersects geometrically once a leaf of $G$ and is otherwise
disjoint from $G$. The following equality  

\begin{equation}
\lbl{eq.join}
[M,\psdraw{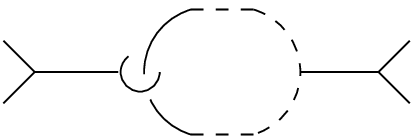}{0.7in}]=[M, \psdraw{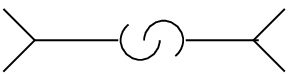}{0.7in}]
+[M,\psdraw{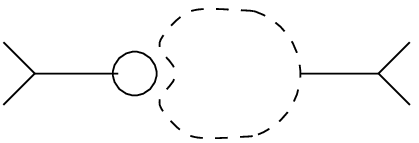}{0.7in}]=[M,\psdraw{nullh2}{0.7in}]
=-[M,\psdraw{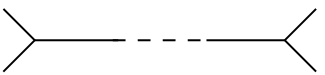}{0.7in}] 
\end{equation}

\noindent
which  follows by Lemma \ref{lem.cut}, concludes our proof.
\end{proof}

\begin{lemma}
\lbl{lem.cbra}
Let $(G,\gamma)$ 
be a \clover \ in a manifold $M$ together with a distinguished leaf $\gamma$
that bounds two surfaces $\S_0$ and $\S_1$ in $M$.
Then,
$$\la G, \S_0-\S_1 \ra=0 \in \GY {}.$$ 
\end{lemma}

\begin{proof}
It follows from two applications of the $\OBRA$ relation that
$$-G=\la G, \S_0 \ra=\la G, \S_1 \ra \in \GY {}.$$
\end{proof}

\begin{proof}(of Theorem \ref{thm.1})
First we construct the map $W_{M,b}$. Let $b=(b_1, \dots, b_r)$ be 
a \spanning \ link in $M$.
Given a graph $G$ with colored legs  choose an arbitrary 
embedding of it in $M$.
For every coloring $\sum_i a_i b_i$ of each of its univalent vertices 
(where $a_i$ are integers), push $|a_i|$ disjoint
copies of $b_i$ (using the framing of $b_i$), orient them the same (resp.
opposite) way from $b_i$ if $a_i \geq 0$ (resp. $a_i < 0$), and finally
take an arbitrary band sum of them. We can arrange the resulting knots,
one for each univalent vertex of the embedding of $G$, to be disjoint
from each other, and together with the embedding of $G$ to form an
embedded graph
with leaves in $M$. Although the isotopy class of the embedded graph
with leaves is not unique, the image of $[M,G] \in \GY {m}$ is well-defined.
This follows from Lemma \ref{lem.cut}.
We need to show that the relations $\AS$, $\IHX$, $\OBRA$, $\CBRA$ 
and $\LOOP$ are
mapped to zero, which will define our map $W_{M,b}$. For the $\AS$ and $\IHX$
relations, see for example \cite[Section 4.1]{GGP}.
The $\LOOP$ relation follows from Figure \ref{loop}.

\begin{figure}[htpb]
$$ \psdraw{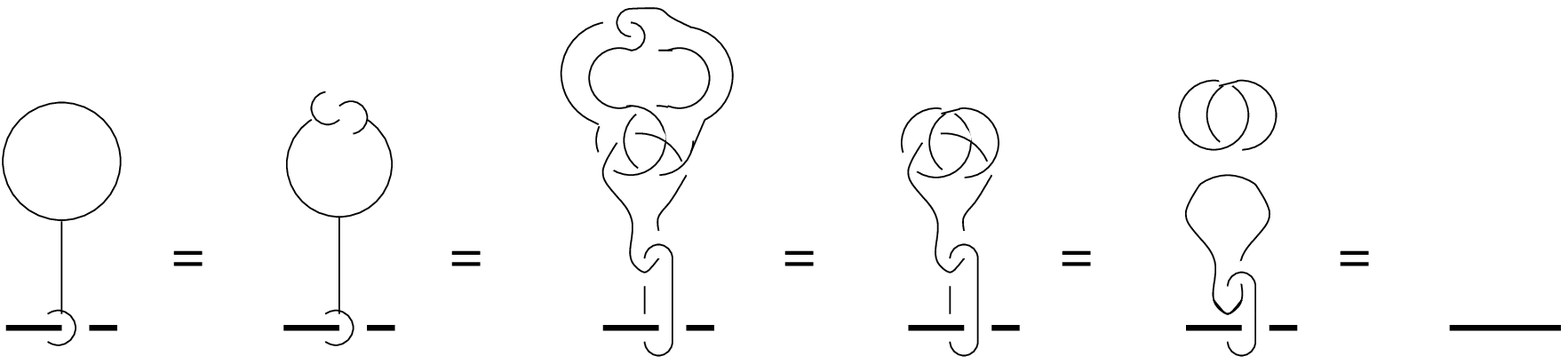}{5.2in} $$
%\caption{.}
\lbl{loop}
\end{figure}

The $\OBRA$ and $\CBRA$ relations follow from Lemmas \ref{lem.obra} and
\ref{lem.cbra}.

We now show that $W_{M,b}$ is onto, over $\BZ[1/2]$. Note first
that $\GY m$ is generated by $[M,G]$ for all simple graphs of degree
$m$, where a simple graph is a disjoint union of graphs of
degree $1$. Each of the leaves of $G$ are isotopic to some connected sum
of (possibly orientation reversed) components of $b$ and contractible
knots. Using the Cutting Lemma \ref{lem.cut}, we may assume that each leaf is
isotopic to one of the components of $b$ (with possibly reversed orientation)
or is contractible in $M$. From this point on, the proof is analogous
to the case of $M=S^3$. Let $L$ be the link consisting of all contractible
leaves of $G$. There exists a trivial, unit-framed link $C$ in $M\smallsetminus
(G \smallsetminus L)$ with the properties that \newline
$\bullet$ each component of $C$ bounds a disk that intersects $L$ at at most
two points. \newline
$\bullet$ Under the diffeomorphism of $M$ with $M_C$, $L$ becomes a 
zero-framed unlink bounding a disjoint collection of disks $D_i$. \newline
Such a link $C$ was called $L$-{\em untying} in \cite{GGP}.
Lemma \ref{lem.cut} and Equations \eqref{eq.join} 
imply that we can assume each of the disks $D_i$ are disjoint
from $G$ and intersect $C$ in at most two points $C$. 
Lemma \ref{lem.relations} imply that $W_{M,b}$ is onto,
over $\BZ[1/2]$.

In order to show that $\Ao (b)$ is independent of $b$, up to isomorphism, we
need the following:

\begin{lemma}
\lbl{lem.1}
Every two \spanning \ links $b$ and $b'$ in $M$ are equivalent by a sequence
of moves 
\begin{itemize}
\item[M1:]
Add one component (after possibly changing its orientation) 
of $b$ to another. 
\item[M2:]
Change the framing of a component of $b$.
\item[M3:]
Insert or delete a null-homologous zero-framed component of $b$.
\end{itemize}
\end{lemma}

\begin{proof}
It suffices to show that under these moves $b$ is equivalent to $b \cup b'$.
Consider a component $b_i'$ of $b'$. Since $b$ is a basis of $H_1(M,\BZ)$,
we can add a multiple of components of $b$ (after perhaps changing their orientation) so that $b_i'$ is nullhomologous,
in which case we can change its framing to zero, and  erase it. 
The lemma now follows by induction on the 
number of components of $b'$.
\end{proof}

If $b'$ is obtained from $b$ by applying one of the three moves above,
we will now define $W_{M,b,b'}:\Ao (b)\to\Ao (b')$ (abbreviated
 by $W_{b,b'}$ in what follows) and show that
Equation \eqref{eq.1} holds.

For the first move, if $b=(b_1,b_2,\dots,b_r)$ and 
$b'=(b_1 \sharp b_2,b_2,\dots,b_r)$ (where $b_1 \sharp b_2$ is 
an arbitrary oriented band sum of $b_1$ with $b_2$)
then $W_{b,b'}$ sends a $b_1$ colored vertex of an abstract graph $G$
to a $b_1 \sharp b_2-b_2$ colored vertex of $G$. It is easy to see
that this defines a map $\Ao (b) \to \Ao (b')$ whose inverse sends a 
$b_1 \sharp b_2$
colored vertex of $G$ to a $b_1+b_2$ colored vertex of $G$.
Similarly, one can define a map $W_{b,b''}$ where
$b'=(b_1 \sharp \overline{b}_2,b_2,\dots,b_r)$.
Equation \eqref{eq.1} follows from Lemma \ref{lem.cut}.

For the second move, let $b=(b_1,b_2,\dots,b_r)$ and
$b'=(b_1',b_2, \dots, b_r)$ where $b_1'$ is a knot whose framing differs
from that of $b_1$ by $\e=\pm 1$. For 
graph $G$ with  $n$ legs colored by $b_1$ we define
$$W_{b,b'}(G)=\sum_{I:|I|=\text{even}}   \e^{|I|/2} G_I'$$
where the summation is over all functions $I:\{1,\dots,n\} \to \{0,1\}$
such that the cardinality $I$ of $I^{-1}(1)$ is even
and $G_I'$ is the result of gluing the $b_1$ colored
legs $l_i$ of $G$ for which $I(i)=1$
pairwise and recoloring the remaining $b_1$ colored legs with
$b_1'$ colored legs. It is easy to see that $W_{b,b'}$ is well-defined (i.e.,
that it respects the relations in $\Ao (b)$) and that its inverse
is given by 
$$W_{b,b'}(G')=\sum_{I:|I|=\text{even}}  (- \e)^{|I|/2} G_I.$$

Let $C$ denotes a $(-\e)$-framed unknot in $M$ which bounds a disk
that geometrically intersects $b_1$ in one point and intersects no other
components of $b$. Then $M_C$ is diffeomorphic to $M$ under a diffeomorphism
that sends the image of $b$ in $M_C$ to $b'$ in $M$.
Since $W_{M,b'}(G)=[M_C,G]$ and $W_{M,b}(G_I)=[M,G_I]$, 
Equation \eqref{eq.1} (or rather, its equivalent form $W_{b'}=
W_{b}\circ W_{b',b}$) follows from the following:

\begin{lemma}
For a graph $G$ of degree $m$ as above, we have
in $\GY {m}$:
$$[M_C,G]=\sum_{I:|I|=\mathrm{even}} (- \e)^{|I|/2} [M,G_I].$$
\end{lemma}

\begin{proof}
Using the Cutting Lemma
\ref{lem.cut} each $b_1$-colored leaf $l_i$ of $G$ can be split along an arc
in two leaves; one that bounds a disk $D_i$ intersecting $C$ once and disjoint
from $b$, and another that is isotopic to $b_1$ but disjoint from $C$.
For $I: \{1, \dots, n\}\to \{0,1\}$, let $G'_I$ denote the graph
$(G \smallsetminus (b_1 \text{ colored leaves of } G))\cup \cup_{i:I(i)=1}
D_i$. Lemma \ref{lem.cut} implies that
$[M_C,G]=\sum_I [M_C,G_I']$. 
Let $G_I''$ denote the graph in $M$ that corresponds to
$G_I'$ under the diffeomorphism $M=M_C$; we obviously have 
$[M_C,G_I']=[M,G_I'']$. 
 Note that $G_I''$ has a collection of $|I|$ 
leaves each of which is unknotted bounding a disk with linking number $\e$ with
every other leaf of this collection. An application of Lemma
\ref{lem.obra} $|I|$ times together with 
Lemma \ref{lem.relations} implies that 
$[M_C,G_I'']=(-\e)^{|I|/2}[M,G_I]$ 
(resp. $0$) for even (resp. odd) $|I|$. 
\end{proof}

For the third move, let $b=(b_1,b_2,\dots,b_r)$ and
$b'=(b_0,b_1,b_2,\dots,b_r)$ where $b_0$ is a null-homologous zero-framed
knot, and consider the natural map $W_{b,b'}:\Ao (b)\to \Ao (b')$.
Choose a surface $\S_0$ that $b_0$ bounds.
The $\OBRA$ relation in $\Ao (b')$ for $b_0$ colored vertices defines
a map $W_{b',b}: \Ao(b')\to\Ao(b)$; this map is independent of $\S_0$
since the difference between two choices of $\S_0$ equals to a choice
of a closed surface and the resulting difference vanishes due to the 
$\CBRA$ relation on $\Ao(b)$. It is easy to see that $W_{b',b}$
 is inverse to $W_{b,b'}$.
Equation \eqref{eq.1} follows essentially by definition.
This completes the proof of Theorem \ref{thm.1}. 
\end{proof}

\begin{proof}(of Corollary \ref{cor.1})
The first statement follows immediately from the fact that if $b$ is a basis
then no nontrivial linear combination is nullhomologous, thus the $\OBRA$
relation is vacuous.

For the second statement, since we are using $\BQ$ coefficients, we may assume
that the link $b$ is a basis for $H_1(M,\BZ)/(\mathrm{torsion})$, and choose
a link $b^t$ to span the torsion part of $H_1(M,\BZ)$.
Then, we have that $\Ac(b)=\Ao(b)\to\Ao(b\cup b^t)$.
There are integers $n_i$ and surfaces $\S_i$ such that $n_i b^t_i=\partial\S_i$
for all components of $b^t$. The $\OBRA$ relation for $b^t$ colored legs
gives a map $\Ao(b\cup b^t)\to\Ao(b)$ which is independent of the choices of
$\{n_i,\S_i\}$ and is inverse to the map $\Ao(b)\to\Ao(b\cup b^t)$.
Thus, $\Ac(b)\cong_{\BQ}\Ao(b\cup b^t)$. Since $\Ao(b\cup b^t)\cong\Ao(b')$
for every \spanning \ link $b'$, the result follows.

The third statement follows immediately from the fact that if the intersection
form on $M$ vanishes, then the $\CBRA$ relation is vacuous.

The forth and fifth statements are immediate consequences of those above.
\end{proof}

\begin{proof}(of Corollary \ref{cor.Q})
Let $b'$ be a \spanning \ link and $b$ be a \qbasis .
Then, we have over $\BQ$
$$
\Ac(H(M))\cong_{\BQ}\Ac(b) \cong_{\BQ} \Ao(b') \to \GY {}$$
which concludes the proof of the corollary.
\end{proof}

\begin{proof}(of Theorem \ref{thm.2})
The proof is a simple application of the {\em locality property}
of the Kontsevich integral, as 
explained leisurely in \cite[II, Section 4.2]{A}, 
and a simple counting argument. 

We now give the details. We need to show that 
\begin{itemize}
\item
The part of the
LMO=Aarhus integral $Z \in \A(\phi)$ of degree
at most $n$ is an invariant of type $n$.
\item
For a trivalent graph $G$ of degree $n$ in a \qhs \ $M$, we have that
$$Z(M_G)=G + \text{ higher degree diagrams } \in \A(\phi).$$
\end{itemize}

For the first claim,
recall that a degree $1$ \clover \ $G$ in a manifold $M$ is the image of
an embedding $ V \to M$ of a neighborhood $V$ of the standard (framed)
graph $\Gamma$ of $\BR^3$, 
and that surgery of $M$ along $G$ can be described as
the result of Dehn surgery on the six component link $L$ in $V$
shown below 
$$\psdraw{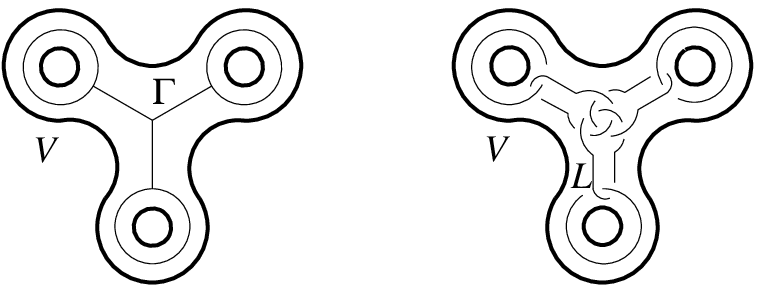}{3.0in}$$
$L$ is partitioned in three blocks $L_1,L_2,L_3$ of two component links
each. We call each block an {\em arm} of $G$.
Alternating a \qhs \ $M$ with respect to
surgery on $G$ equals to alternating $M$ with respect to all nine
subsets of the set of arms of $G$.

Recall also that the Kontsevich integral of a framed link $L$
in a 3-manifold $M$
$Z(M,L)$ (defined by
Kontsevich for links in $S^3$ and extended by Le-Murakami-Ohtsuki
for links in arbitrary 3-manifolds \cite[Section 6.2]{LMO})
takes values in linear combinations of $L$-colored uni-trivalent
graphs.

Recall also that the LMO=Aarhus integral of a \qhs \ $M_L$ (obtained
by surgery on a framed link $L$
in a \qhs \ $M$) is obtained
by considering the Kontsevich integral $Z(M,L)$, 
splitting it in a
quadratic $Z^q$ and trivalent (a better name would be ``other'') 
part $Z^t$, and gluing the $L$-colored legs
of $Z^t$ using the inverse linking matrix of $L$.

Given a \clover \ $G=\cup_{i=1}^{n} G_i$ in a \qhs \ $M$,
(where $G_i$ are of degree $1$), let $L^{\mathrm{act}}$ denote the
link that consists of the $3n$ arms of $G$. When we compute $Z([M,G])=
Z([M,L^{\mathrm{act}}])$, we need to concentrate on all 
$L^{\mathrm{act}}$-colored uni-trivalent graphs 
that have at least one univalent vertex on each block of $G$.
Such graphs will have at least $3n$ univalent vertices.
Since at most three univalent vertices can share a trivalent vertex, 
it follows that the above considered graphs will have at least $n$
trivalent vertices; in other words it follows that
$Z([M,G]) \in \A_{\geq n}(\phi)$.

The second claim is best shown by example. Recall that surgery on the
(generic trivalent graph) 
$\Theta$ shown 
below corresponds to surgery on two \clover s $G_1$ and $G_2$, each with
arms $\{E_{ij},L_{ij}\}$ for $i=1,2$ and $j=1,2,3$.
The linking matrix of the $12$ component link $L^{\mathrm{act}}=
E_{ij} \cup L_{ij}$
and its inverse are given by
$$
\left(\begin{array}{cc}
0 & I \\ I & I \\
\end{array}\right)
\text{ and }
\left(\begin{array}{cc}
-I & I \\ I & 0 \\
\end{array}\right)
$$
where $I$ is the identity $6\times 6$ matrix. 
The relevant trivalent part $Z^t(M,L^{\mathrm{act}})$ 
is shown schematically in four cases here, where the graphs on the left terms
of each case 
come from $G_1$ and the graphs on the right terms of each case come from $G_2$
and the dashed lines correspond to gluings of the univalent vertices:
$$\psdraw{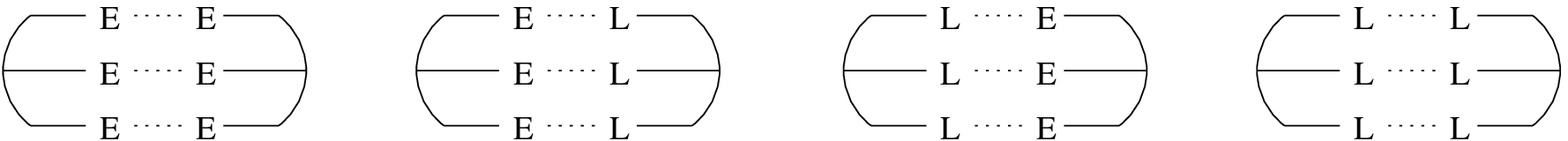}{5in}$$
However, the last three cases all contribute zero, since $LLL$ is
a 3-component unlink whose coefficient in $Z^t$ is a multiple of
the triple Milnor invariant and thus vanishes.
Thus, we are only left to glue terms in the first case, and this is
summarized in the following figure 
$$\psdraw{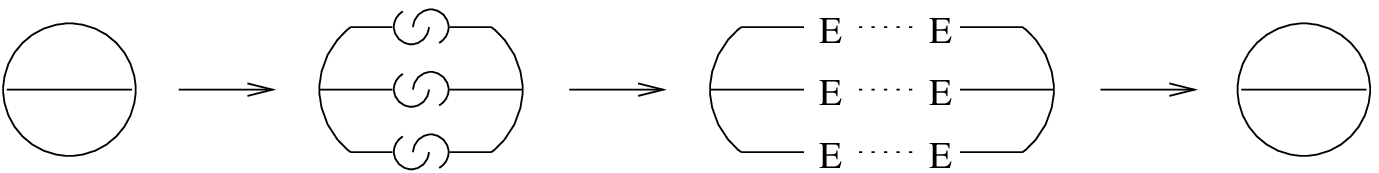}{3in}$$
which concludes the proof.
\end{proof}

\begin{proof}(of Corollary \ref{cor.3})
If $G$ is as in the statement of the corollary, colored by a sublink $b'$
of $b$ which is not \spanning \ , then we can find an $x \in b\smallsetminus 
b'$, and a closed surface $x^\ast$ such that
$[x^\ast][y]=\delta_{y,x}$ for all components $y$ of $b$. 
Cut $G$ along an edge, and color the two new leaves $x$ and $\ast$ to
obtain a graph $(G,\ast)$. By definition, we have $\langle G, x^\ast \rangle
= G$, thus the result follows from the $\CBRA$ relation.
\end{proof}

\begin{proof}(of Corollary \ref{cor.levine})
We will first show the result for $k=0$.
Let $G$ be as in the statement of the corollary and let $G'$
be the graph with two more leaves than $G$, colored by $x$ and $\ast$
respectively as shown: 
$$\psdraw{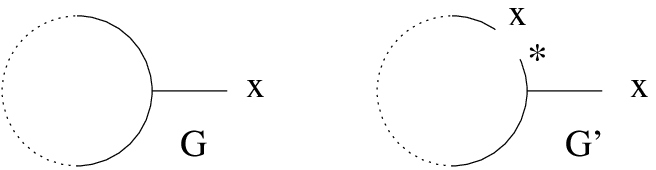}{2.5in}.$$
The $\CBRA$ relation implies that 
$$0=\langle G',x^\ast
\rangle = G + \LOOP = G=G^{(0)}.$$
Now, we will show the result for all $k$.
Let $G(n)$ be the same graph as $G$ with $r+1$ leaves colored by 
$x, x+ny_1,\dots,x+ny_r$, for $n \in \BN$. Since $x$ is primitive and
linearly independent from $\{x+ny_1,\dots,x+ny_r\}$, the $k=0$ case
for $G(n)$ shown above implies
that $G(n)=0$ for all $n$. Since $G(n)=\sum_k n^k G^{(k)}$, the result follows.
\end{proof}

\noindent
{\bf Acknowledgement}
We wish to thank D. Bar-Natan, J. Levine, M. Polyak and S. Weinberger 
for stimulating conversations.

\ifx\undefined\bysame
	\newcommand{\bysame}{\leavevmode\hbox
to3em{\hrulefill}\,}
\fi

\end{document}